\documentclass[10pt,letterpaper]{article}
\makeatletter
\usepackage{amsfonts}
\usepackage{amsmath,amsthm,amssymb}

\newtheorem{thm}{\textbf Theorem}[section]

\newtheorem{rem}{\textbf Remark}[section]
\newtheorem{cor}{\textbf Corollary}[section]

\newtheorem{prop}{\textbf Proposition}[section]
\newtheorem{defin}{\textbf Definition}[section]

\newcommand{\md}{\mbox{d}}
\newcommand{\be}{\begin{eqnarray}}
\newcommand{\ee}{\end{eqnarray}}
\newcommand{\mr}{\mathbb{R}}

\newcommand{\mx}{\mbox}

\newcommand{\pt}{\partial}
\newcommand{\bes}{\begin{eqnarray*}}
\newcommand{\ees}{\end{eqnarray*}}

\newcommand{\om}{\omega}

\newcommand{\ka}{\kappa}

\newcommand{\De}{\Delta}
\newcommand{\ga}{\gamma}

\begin{document}
\begin{titlepage}
\title{\bf On the regularity criteria of weak solutions to the micropolar fluid equations in Lorentz space}
\author{ Baoquan Yuan
       \\ School of Mathematics and Informatics,
       \\ Henan Polytechnic University, Henan, 454000, China.\\
        (bqyuan@hpu.edu.cn)
          }
\date{}
\end{titlepage}
\maketitle

%

%
%

\begin{abstract}
In this paper the regularity of weak solutions and the blow-up
criteria of smooth solutions to the micropolar fluid equations on
three dimension space are studied in the Lorentz space
$L^{p,\infty}(\mathbb{R}^3)$. We obtain that if $u\in
L^q(0,T;L^{p,\infty}(\mathbb{R}^3))$ for $\frac2q+\frac3p\le 1$ with
$3<p\le \infty$; or $\nabla u\in
L^q(0,T;L^{p,\infty}(\mathbb{R}^3))$ for $\frac2q+\frac3p\le 2$ with
$\frac32<p\le \infty$; or the pressure $P\in
L^q(0,T;L^{p,\infty}(\mathbb{R}^3))$ for $\frac2q+\frac3p\le 2$ with
$\frac32<p\le \infty$; or $\nabla P\in
L^q(0,T;L^{p,\infty}(\mathbb{R}^3))$ for $\frac2q+\frac3p\le 3$ with
$1<p\le \infty$, then the weak solution $(u,\omega)$ satisfying the
energy inequality is a smooth solution on $[0,T)$.
 \vskip0.1in


\end{abstract}

\vspace{.2in} {\bf Key words:}\quad Micropolar fluid equations,
regularity of weak solutions, Lorentz spaces.


\section{Introduction}
\setcounter{equation}{0}

This paper concerns about the regularity of weak solutions and
blow-up criteria of smooth solutions to the micropolar fluid
equations in three dimensions
 \be\label{MMP}
 \begin{cases}
\frac{\pt u}{\pt t}-(\mu+\chi)\De u+u\cdot\nabla u+\nabla P-\chi\nabla\times\omega=0,\\
\frac{\pt\omega}{\pt t}-\gamma\Delta \om-\kappa \nabla div
\omega+2\chi \omega+u\cdot\nabla\omega-\chi\nabla\times u=0,\\
div u=0,\\
u(x,0)=u_0(x),\ \omega(x,0)=\omega_0(x),
 \end{cases}
 \ee
 where $u=(u_1(t,x),u_2(t,x),u_3(t,x))$ denotes the velocity of
the fluids at a point $x\in \mr^3$, $t\in [0,T)$,
$\om=(\om_1(t,x),\om_2(t,x),\om_3(t,x))$ and $P=P(t,x)$ denote,
respectively, the micro-rotational velocity and the hydrostatic
pressure. $u_0$, $\om_0$ are the prescribed initial data for the
velocity and angular velocity with properties div $u_0=0$. $\mu$
is the kinematic viscosity, $\chi$ is the vortex viscosity,
$\kappa$ and $\gamma$ are spin viscosities. Theory of micropolar
fluids was first proposed by Eringen \cite{Eringen} in 1966, which
enables us to consider some physical phenomena that cannot be
treated by the classical Navier-Stokes equations for the viscous
incompressible fluids, for example, the motion of animal blood,
liquid crystals and dilute aqueous polymer solutions etc. The
existences of weak and strong solutions were treated by Galdi and
Rionero \cite{G-R}, and Yamaguchi \cite{Yamaguchi}, respectively.
If, further, the vortex viscosity $\chi=0$, the velocity $u$ does
not depend on the micro-rotation field $\om$, and the first
equation reduces to the classical Navier-Stokes equation which has
been greatly analyzed, see, for example, the classical books by
Ladyzhenskaya \cite{Ladyzhenskaya}, Lions \cite{Lions} or
Lemari\'{e}-Rieusset \cite{Lemarie}.

There is a large number of literature on the mathematical theory
of micropolar fluid equations (\ref{MMP}) (see, for example,
\cite{Lange,Yamaguchi,G-R,F-V,D-C,V-R,C-D-C,C-C-D}). The existence
and uniqueness of global solutions were extensively studied by
Lange \cite{Lange}, Galdi and Rionero \cite{G-R}, Yamaguchi
\cite{Yamaguchi}. Recently, Ferreira and Villamizar-Roa \cite{F-V}
considered the existence and stability of solutions to the
micropolar fluids in exterior domains. Villamizar-Roa and
Rodr\'{\i}guez-Bellido \cite{V-R} studied the micropolar system in
a bounded domain using the semigroup approach in $L^p$, showing
the global existence of strong solutions for small data and the
asymptotic behavior and stability of the solutions. Concerning the
dynamic behavior of solutions to equations (\ref{MMP}) one may
refer to the references \cite{C-D-C,C-C-D,D-C} and references
therein.

The purpose of this paper is to study the regularity of weak
solutions and the breakdown criteria of smooth solutions to the
micropolar fluid equations (\ref{MMP}). The classical blow-up
criteria of smooth solutions to the Navier-Stokes equations also
hold for the micropolar fluid equations. For the Navier-Stokes
equations, Serrin \cite{Serrin}, prodi \cite{Prodi} and Veiga
\cite{Veiga} established the classic Serrin-type regularity criteria
of weak solutions in terms of $u$ or its gradient $\nabla u$. Later,
many improvements and extensions were established, for example, see
\cite{K-T1,K-O-T,K-S,Y-Z} and references therein. Berselli and Galdi
\cite{B-G}, Chae and Lee \cite{C-L} obtained the regularity criteria
of weak solutions in terms of the pressure $P$ or its gradient
$\nabla P$. Later, Zhou Y improved it in terms of the pressure and
its gradient in a general domain \cite{Zhou1}, and Zhou Y
\cite{Zhou2,Zhou3}, Struwe \cite{Struwe} obtained the regularity
criteria of weak solutions in terms of the gradient of pressure.

\begin{thm}
Suppose $u\in L^\infty([0,T]);L^2(\mr^n))\cap
L^2([0,T];H^1(\mr^n))$ is a Leray-Hopf weak solution to the
Navier-Stokes equations, $P$ is the pressue. If one of the
following conditions is satisfied:

{\rm (1)} $u\in L^q(0,T;L^p(\mr^n))$ for $\frac2q+\frac np\le 1$
with $n< p\le\infty$;

{\rm (2)} $\nabla u\in L^q(0,T;L^p(\mr^n))$ for $\frac2q+\frac np\le
2$ with $\frac n2<p\le\infty$;

{\rm (3)} $P\in L^q(0,T;L^p(\mr^n))$ for $\frac2q+\frac np\le 2$
with $\frac n2<p\le\infty$;

{\rm (4)} $\nabla P\in L^q(0,T;L^p(\mr^n))$ for $\frac2q+\frac
np\le 3$ with $\frac n3<p<\infty$.

Then $u$ is a smooth solution on $[0,T)$.
\end{thm}

Similar to the Serrin type regularity criteria, Yuan \cite{Yuan}
established the regularity criteria of weak solutions to
magneto-micropolar equations, which is the microploar equations
(\ref{MMP}) coupled with magnetic field $b$, as follows.

\begin{thm}\label{thm0.1}
Let $(u_0,\om_0,b_0)\in L^2(\mr^3)$ with $div u_0=div b_0=0$.
Assume that $(u,\om,b)$ is a Leray-Hopf type weak solution to the
magneto-micropolar equations. If one of the following conditions
holds:

{\rm (1)} $u\in L^q(0,T;L^p(\mr^3))$ for $\frac2q+\frac3p\le 1$
with $3<p\le\infty$;

{\rm (2)} $\nabla u\in L^q(0,T;L^p(\mr^3))$ for $\frac2q+\frac3p\le
2$ with $\frac32<p\le\infty$.

Then $(u,\om,b)$ is a smooth solution on $[0,T)$ with the initial
value $(u_0,\om_0,b_0)$.
\end{thm}

It is worthy to note that the regularity conditions of weak
solutions to magneto-micropolar equations
 are only imposed on the velocity field $u$, which is very important.
 For the magneto-hydrodynamic equations, He and Xin \cite{H-X} first studied and
 established the regularity criteria only imposed on the velocity field
 $u$ or its gradient $\nabla u$. Later, Zhou Y \cite{Zhou4} improved the regularity criteria imposed only  on
 $u$ or its gradient $\nabla u$; and He-Wang \cite{H-W} improved it to the
 weak $L^p$ spaces imposed only on $u$ or its gradient $\nabla u$; and Chen-Miao-Zhang \cite {C-M-Z}
 also improved it to the more general Besov-type space on
 Littlewood-Paley decomposition imposed only on $\nabla\times u$.
 The regularity criteria
of weak solutions to the system (\ref{MMP}) play a important role to
understanding the physical essence of the micropolar fluid motion.
The aim of this paper is to prove that to secure the regularity of
weak solutions to (\ref{MMP}), one only needs to impose conditions
on the velocity field $u$ or its gradient $\nabla u$ or the pressure
of the fluids in the Lorentz spaces. In details, one only need one
of the following conditions to prove the regularity of weak
solutions $(u,\om)$ on $[0,T]$.

{\rm (1)} $u\in L^q(0,T;L^{p,\infty}(\mr^3))$ for $\frac2q+\frac
3p\le 1$ with $3< p\le\infty$;

{\rm (2)} $\nabla u\in L^q(0,T;L^{p,\infty}(\mr^3))$ for
$\frac2q+\frac3p\le 2$ with $\frac32<p\le\infty$;

{\rm (3)} $P\in L^q(0,T;L^{p,\infty}(\mr^3))$ for $\frac2q+\frac
3p\le 2$ with $\frac 32<p\le\infty$;

{\rm (4)} $\nabla P\in L^q(0,T;L^{p,\infty}(\mr^3))$ for
$\frac2q+\frac 3p\le 3$ with $1<p<\infty$.

This demonstrates that, in the regularity of weak solutions, the
micro-rotational velocity $\om$ of particles play less important
role than the velocity $u$ does, and the regularity of weak
solutions to (\ref{MMP}) is dominated by the velocity $u$ of the
fluids.

We conclude this introduction by describing the plan of the paper.
We give our main results of the blow-up criteria for a smooth
solution to (\ref{MMP}) and as applications we prove the
regularity of weak solutions in section 2. Section 3 devoted to
prove Theorem \ref{thm2.1} and \ref{thm2.2}, respectively.

\section{Main results}
\setcounter{equation}{0}
 Before stating our main results we introduce some function spaces,
 notations and generalized H\"older inequality.
 Let
$C^\infty_{0,\:\sigma}(\mathbb{R}^3)$ denote the set of all
$C^\infty$ vector functions $f(x)=(f_1(x),f_2(x),f_3(x))$ with
compact support such that div $ f(x)=0$.
$L^r_\sigma(\mathbb{R}^3)$ is the closure of
$C^\infty_{0,\:\sigma}(\mathbb{R}^3)$-function with respect to the
$L^r$-norm $\|\cdot\|_r$ for $1\le r\le\infty$.
$H^s_\sigma(\mathbb{R}^3)$ denotes the closure of
$C^\infty_{0,\:\sigma}(\mathbb{R}^3)$ with respect to the
$H^s$-norm $\|f\|_{H^s}=\|(1-\De)^{\frac s2}f\|_2$, for $s\ge 0$.

In the following part we recall Lorentz spaces. Let
$(X,\mathcal{M},\mu)$ be a non-atomic measurable space. For the
complex-valued or real-valued, $\mu$-measurable function $f(x)$
defined on $X$, its distribution function is defined by
\label{fm3}\begin{eqnarray} f_{*}(\sigma)=\mu \{x\in X :|f(x)|>
\sigma \},\ for\ \sigma>0,
\end{eqnarray}
which is non-increasing and continuous from the right.
Furthermore, its non-increasing rearrangement $f^{*}$ is defined
by \label{nir}\begin{eqnarray} f^{*}(t)=\inf\{s>0:f_{*}(s)\leq
t\},\ for\ t>0.
\end{eqnarray}
which is also non-increasing and continuous from the right, and has
the same distribution function as $f(x)$.

The Lorentz space $L^{p,q}$ on $(X,\mathcal{M},\mu)$ is the
collection of all the real-valued or complex-valued,
$\mu-$measurable function $f(x)$ defined on $X$ such that
$\|f\|_{p,q}<\infty$ with
\begin{eqnarray}
 \|f\|_{p,q}=\begin{cases}
   \bigg(\frac{q}{p}\displaystyle\int^\infty_0
   (t^{\frac{1}{p}}f^*(t)\big)^q\frac{\mathrm{d}t}{t}\bigg)^{\frac{1}{q}},
   &\ if\ 1\leq p<\infty,\  1<q<\infty,\\
  \sup\limits_{\:t>0}t^{\frac{1}{p}}f^*(t),
  &\ if\ 1\leq p\leq\infty,\ q=\infty.\end{cases}
\end{eqnarray}
If $q=\infty$ write $L^{p,\infty}(\mr^3)$ as $L^p_w(\mr^3)$ which
is the weak $L^p$ space. Moreover
 \begin{eqnarray}\label{2.14}
\|f\|_{p,\infty}=\sup\limits_{\:t>0}t^{\frac{1}{p}}f^{*}(t)=
\sup\limits_{\:\alpha>0}\alpha (f_{*}(\alpha))^{\frac{1}{p}},
\end{eqnarray}
for any $f(x)\in L^{p,\infty}$. For details see \cite{L,M1} and
\cite{S-W}.

We also need the H\"older inequality in Lorentz spaces which we
recall as follows, for details see O'Neil \cite{Oneil}.

\begin{prop}
Let $1< p_1,\ p_2,\ r<\infty$ satisfying
 \begin{eqnarray*}
\frac1{p_1}+\frac1{p_2}<1,\ \frac1r=\frac 1{p_1}+\frac 1{p_2},
 \end{eqnarray*}
 and $1\le q_1,\ q_2,\ s\le\infty$ satisfying
\begin{eqnarray*}
\frac1{q_1}+\frac1{q_2}\ge\frac1s.
\end{eqnarray*}
If $f\in L^{p_1, q_1}$ and $g\in L^{p_2,\:q_2}$, then $fg\in L^{r,
s}$, and the generalized H\"older inequality
\begin{eqnarray}\label{holder}
\|h\|_{r,s}\le r'\|f\|_{p_1,q_1}\|g\|_{p_2,q_2}
\end{eqnarray}\\
holds, where $r'$ stands for the dual to $r$, i.e.
$\frac1r+\frac1{r'}=1$.
\end{prop}

To this end, we state the main results as follows.

\begin{thm}\label{thm2.1}
Let $u_0(x)\in H^1_\sigma(\mr^3)$ and $\om_0(x)\in H^1(\mr^3)$.
Assume that $u(t,x)\in C([0,T);H^1_\sigma(\mr^3))\cap
C((0,T);H^2_\sigma(\mr^))$ and $\om(t,x)\in C([0,T);H^1(\mr^3))\cap
C((0,T);H^2(\mr^))$ is a smooth solution to the
equations(\ref{MMP}). If $u$ satisfies one of the following
conditions

{\rm (a)} $u(t,x)\in L^q((0,T);L^{p,\infty}(\mr^3)),\ \mx{ for }
\frac2q+\frac3p\le 1\ \mx{ with }3<p\le\infty$;

{\rm (b)} $\nabla u(t,x)\in L^q((0,T);L^{p,\infty}(\mr^3)),\ \mx{
for } \frac2q+\frac3p\le 2\ \mx{ with }\frac32<p\le\infty$.

Then the solution $(u,\om)$ can be extended smoothly to $[0,T')$ for
some a $T'>T$.
\end{thm}


\begin{thm}\label{thm2.2}
Let $u_0(x)\in L^4_\sigma(\mr^3)$ and $\om_0(x)\in L^4(\mr^3)$.
Assume that $u(t,x)\in C([0,T);L^4_\sigma(\mr^3))\cap
C((0,T);H^{1,4}_\sigma(\mr^3))$ and $\om(t,x)\in
C([0,T);L^4(\mr^3))\cap C((0,T);H^{1,4}(\mr^3))$ is a smooth
solution to the equations (\ref{MMP}), and $P$ is the pressure. If
$P$ satisfies the condition

{\rm (1)} $P(t,x)\in L^q((0,T);L^{p,\infty}(\mr^3)),\ \mx{ for }
\frac2q+\frac3p\le 2\ \mx{ with } \frac32<p\le\infty$;

Or the gradient of the pressure $\nabla P$ satisfies the condition

{\rm (2)} $\nabla P(t,x)\in L^q((0,T);L^{p,\infty}(\mr^3)),\ \mx{
for } \frac2q+\frac3p\le 3\ \mx{ with } 1<p\le\infty$.

Then the solution $(u,\om)$ can be extended smoothly beyond $t=T$.
\end{thm}

We next consider the criteria of regularity of weak solutions to
the micropolar equations (\ref{MMP}), for this purpose we first
introduce the definition of a weak solution.

\begin{defin}
Let $u_0(x)\in L^2_\sigma(\mathbb{R}^3)$ and $\om_0(x)\in
L^2(\mr^3)$. A measurable function $(u(t,x),\om(t,x))$ is called a
weak solution to the micropolar equations (\ref{MMP}) on $[0,T)$,
if
\\(a)
\bes u(t,x)\in L^\infty([0,T);L^2_\sigma(\mr^3))\cap
L^2([0,T);H^1_\sigma(\mr^3)),
 \ees
 and
 \bes
 \om(t,x)\in
L^\infty([0,T);L^2(\mr^3))\cap L^2([0,T);H^1(\mr^3));
 \ees
 (b)
 \bes
  &&\int_0^T\{-(u, \partial_{\tau}\varphi)+(\mu+\chi)(\nabla u,\nabla
  \varphi)-
(u\cdot\nabla \varphi, u)+\chi (\nabla\times\varphi,\om) \md\tau
\\&&=(u_0, \varphi(0)),
  \ees
\bes &&\int_0^T\{-(\om,\pt_\tau\phi)+\gamma(\nabla
\om,\nabla\phi)+\kappa(div\om,div\phi)+2\chi(\om,\phi)
-(u\cdot\nabla \phi,\om)+\chi(\nabla\times
\phi,u)\md\tau\\&&=(\om_0,\phi(0)), \ees
  for any $\varphi(t,x)\in
H^1[(0,T);H^1_{\sigma}(\mr^3)$ and $\phi(t,x)\in
H^1([0,T);H^1(\mr^3)$ with $\varphi(T)=0$ and $\psi(T)=0$.
\end{defin}
In the reference \cite{RM-B}, Rojas-Medar and Boldrini proved the
global existence of weak solutions to the equations of the
magneto-micropolar fluid motion  by the Galerkin method. The weak
solutions $(u,\om)$ also satisfy the energy inequality
 \be\label{EI}
&&\|(u,\om)\|^2_2+2\mu\int^t_0\|\nabla u\|^2_2\md
s+2\ga\int^t_0\|\nabla\om\|^2_2\md s\\ \nonumber
&&+2\ka\int^t_0\|\mx{div} \om\|^2_2\md
s+2\chi\int^t_0\|\om\|^2_2\md s\le \|(u(0),\om(0))\|^2_2,
 \ee
for $0< t\le T$.

As immediate corollaries we establish the regularity criteria of
weak solutions.

\begin{cor}\label{cor2.1}
Let $u_0(x)\in H^1_\sigma(\mr^3)$ and $\om_0(x)\in H^1(\mr^3)$.
Assume that $(u(t,x),\om(t,x))$ is a weak solution to the
equations(\ref{MMP}) and satisfies the energy inequality
(\ref{EI}). If $u$ satisfies one of the following conditions

{\rm (a)} $u(t,x)\in L^q((0,T);L^{p,\infty}(\mr^3)) \mx{ for }
\frac2q+\frac3p\le 1\ \mx{ with }3<p\le\infty$;

{\rm (b)} $\nabla u(t,x)\in L^q((0,T);L^{p,\infty}(\mr^3)) \mx{
for } \frac2q+\frac3p\le 2\ \mx{ with }\frac32<p\le\infty$.

Then the solution $(u,\om)$ is a regular solution on $(0,T]$.
\end{cor}


\begin{cor}\label{cor2.2}
Let $u_0(x)\in L^2_\sigma\cap L^4_\sigma(\mr^3)$ and $\om_0(x)\in
L^2\cap L^4(\mr^3)$. Assume that $(u(t,x),\om(t,x))$ is a weak
solution to the equations (\ref{MMP}) and satisfies the energy
inequality (\ref{EI}), and $P$ is the pressure. If $P$ satisfies
the condition

{\rm (1)} $P(t,x)\in L^q((0,T);L^{p,\infty}(\mr^3)),\ \mx{ for }
\frac2q+\frac3p\le 2\ \mx{ with } \frac32<p\le\infty$;

Or the gradient of the pressure $\nabla P$ satisfies the condition

{\rm (2)} $\nabla P(t,x)\in L^q((0,T);L^{p,\infty}(\mr^3)),\ \mx{
for } \frac2q+\frac3p\le 3\ \mx{ with } 1<p\le\infty$.

Then the solution $(u,\om)$ is a regular solution on $(0,T]$.
\end{cor}

\begin{rem}
In Theorems \ref{thm2.1} and \ref{thm2.2}, if $p=\infty$, then the
 space $L^{\infty,\infty}$ is identical with $L^\infty$.
\end{rem}

\begin{rem}
In this paper the regularity Theorem \ref{thm2.1} and Corollary
\ref{cor2.1} are established for the solution of the micropolar
equations (\ref{MMP}). By the coupling of velocity field and
magnetic field, the conclusions of Theorem \ref{thm2.1} and
Corollary \ref{cor2.1} are also  valid in the magneto-micropolar
equations.
\end{rem}

\begin{rem}
For the magneto-hydrodynamic equations, He-Wang \cite{H-W} proved
that to assure the regularity of weak solution one only need the
condition
 \be\label{hw1}
\nabla u\in L^q((0,T);L^{p,\infty}(\mr^3)),
 \ee
for $\frac2q+\frac3p=2$ with $1<q\le 2$. In our results Theorem
\ref{thm2.1} (b) and Corollary \ref{cor2.2} (b) $1\le q<\infty$ is
more general than that in (\ref{hw1}). Moreover, our proof which is
based on the a priori estimate of the $H^1-$norm of solution is
simple, while it was based on the a priori estimate of the
$L^p-$norm for $p\ge 3$.
\end{rem}

The proofs of Corollaries \ref{cor2.1} and \ref{cor2.2} are
standard. For the completeness we sketch out the proof of
Corollary \ref{cor2.1} only. Since $u_0(x)\in H^1_\sigma(\mr^3)$
and $\om_0(x)\in H^1(\mr^3)$, by the local existence theorem of
strong solution to the micropolar equations (\ref{MMP}), there
exists a unique solution $(\hat u,\hat\om)$ satisfying that $\hat
u(t,x)\in C([0,T^*);H^1_\sigma(\mr^3))$ and $\hat \om(t,x)\in
C([0,T^*);H^1(\mr^3))$ on a small time interval $[0,T^*)$. Since
$(u,\om)$ is a weak solution satisfying the energy inequality
(\ref{EI}), it follows from the Serrin type uniqueness criterion
\cite{Serrin} that $(u(t),\om(t))\equiv (\hat u(t), \hat\om(t))$
on $[0,T^*)$. Thus it is sufficient to show that $T=T^*$. If not,
suppose that $T^*<T$. Without loss of generality, one may assume
that $T^*$ is the maximal existence time of the strong solution
$(\hat u,\hat \om)$. By the conditions (a) or (b) in Corollary
\ref{cor2.1}, we have
 \bes
\int^{T^*}_0\|\hat u(t)\|_{p,\infty}^q\md t<\infty, \mx{ for
}\frac2q+\frac3p\le 1\ \mx{ with }3<p\le\infty,
 \ees
or
 \bes
\int^{T^*}_0\|\nabla\hat u(t)\|_{p,\infty}^q\md t<\infty, \mx{ for
}\frac2q+\frac3p\le 2\ \mx{ with }\frac32<p\le\infty
 \ees
because of $(\hat u(t),\hat \om(t))\equiv ((u(t),\om(t)))$.
Therefore it follows from Theorem \ref{thm2.1} that there exists a
time $T'>T^*$ such that $(\hat u,\hat \om)$ can be extended
smoothly to $[0,T')$, which contracts to the maximality of $T^*$.
We thus complete the proof of Corollary \ref{cor2.1}.

In the following arguments the letter $C$ denotes inessential
constants which may vary from line to line, but does not depend on
particular solutions or functions. We also use
$C(\chi,\ga,\cdots)$ to denote a constant which depends on the
parameters $\chi,\ga,\cdots$ and may vary from line to line.

\section{Proof of Theorems \ref{thm2.1} and \ref{thm2.2}}
\setcounter{equation}{0}

In this section we prove Theorem \ref{thm2.1} and \ref{thm2.2} by a
simple method.

 {\bf Proof of Theorem \ref{thm2.1}}: We differentiate the
equations (\ref{MMP}) with respect to $x_i$, then multiply the
resulting equations by $\pt_{x_i}u,\ \pt_{x_i}\om$ for $i=1,\ 2,\
3$, respectively, integrate with respect to $x$ and sum them up. It
follows that
 \be\nonumber
  && \frac12\frac{\md}{\md
t}(\|(\pt_{x_i}u,\pt_{x_i}\om)\|^2_2)
+\sum_{j=1}^3\Big((\mu+\chi)\|\pt^2_{x_ix_j}u\|^2_2+\ga\|\pt^2_{x_ix_j}\om\|^2_2\Big)
+\kappa\|div\pt_{x_i}\om\|^2_2+ 2\chi\|\pt_{x_i}\om\|^2_2 \\
&\le& \nonumber |(\pt_{x_i}u\cdot\nabla
u,\pt_{x_i}u)|+|(\pt_{x_i}u\cdot\nabla \om,\pt_{x_i}\om)|+2\chi|(\nabla\times\pt_{x_i}u,\pt_{x_i}\om)|\\
\label{3.10} &=&I_1+I_2+I_3,
 \ee
 where we have used the facts that
 \be (\nabla\times\pt_{x_i}u,\pt_{x_i}\om)=
(\nabla\times\pt_{x_i}\om,\pt_{x_i}u)
 \ee
and
 \be (u\cdot\nabla \pt_{x_i}u,\pt_{x_i}u)=(u\cdot\nabla
\pt_{x_i}\om,\pt_{x_i}\om)=0,
 \ee
where $(\cdot,\cdot)$ denotes the $L^2$ inner product on $\mr^3$.
 For conciseness, the
short notation
 \be
 \|(A,B)\|^2_2=\|A\|^2_2+\|B\|^2_2
 \ee
has been used and will be used in the following parts.

(a) We estimate the term $I_1,\ I_2$ and $I_3$ respectively. By
integrations by parts and the generalized H\"older inequality
(\ref{holder}) it follows that
 \be\nonumber
 I_1&\le& \Big|\int_{\mr^3}\pt_{x_i}u\cdot\nabla\pt_{x_i}u\cdot u(x)\md
 x\Big|+\Big|\int_{\mr^3}\pt_{x_i}\pt_{x_i}u\cdot\nabla u\cdot u(x)\md
 x\Big|\\ \nonumber &\le& C(p)\|u\|_{p,\infty}\Big(\|\nabla
 \pt_{x_i}u\pt_{x_i}u\|_{\frac{p}{p-1},1}+\|\pt_{x_ix_i}u\pt_{x_i}u\|_{\frac{p}{p-1},1}\Big)
 \\\label{3.11} &\le& C(p)\|u\|_{p,\infty}\|\nabla
 u\|_{\frac{2p}{p-2},2}\|D^2u\|_2.
 \ee
Applying the real interpolation (see \cite{B-L})
 \bes
L^{\frac{2p}{p-2},2}(\mr^3)=(L^2,L^6)_{\frac{p-3}{p},2}(\mr^3)
 \ees
and the Sobolev embedding $L^6(\mr^3)\hookrightarrow \dot
H^1(\mr^3)$ it yields that
 \be\label{3.12}
\|\nabla u\|_{\frac{2p}{p-2},2}\le C\|\nabla
u\|_2^{1-\frac3p}\|D^2u\|_2^{\frac3p}.
 \ee
Inserting the above estimate (\ref{3.12}) into the estimate
(\ref{3.11}) of $I_1$ one has
 \be\nonumber
I_1&\le& C(p)\|u\|_{p,\infty}\|\nabla
u\|_2^{1-\frac3p}\|D^2u\|^{1+\frac3p}\\&\le& \label{3.13}
\frac\chi{12}\|D^2u\|^2+C(p,\chi)\|u\|_{p,\infty}^\frac{2p}{p-3}\|\nabla
u\|_2^2.
 \ee
Similarly, for $I_2$ one also has
 \be\label{3.14}
I_2\le\frac\ga6\|D^2\om\|^2+C(p,\ga)\|u\|_{p,\infty}^\frac{2p}{p-3}\|\nabla
\om\|_2^2.
 \ee
For the term $I_3$, H\"older and Young inequalities imply that
 \be\label{I3}
I_3\le \frac\chi
2\|\nabla\times\pt_{x_i}u\|^2_2+2\chi\|\nabla\om\|^2_2.
 \ee
Inserting the estimates (\ref{3.13})-(\ref{I3}) into the inequality
(\ref{3.10}) and summing up $i$ from $1$ to $3$, it follows that
 \bes
&&\frac{\md}{\md t}(\|(\nabla u,\nabla\om)\|^2_2)
+(2\mu+\frac12\chi)\|D^2u\|^2_2+\ga\|D^2\om\|^2_2+2\kappa\|\nabla
div\om\|^2_2\\ &\le&
C(p,\chi,\ga)\|u\|_{p,\infty}^{2p/(p-3)}\|(\nabla
u,\nabla\om)\|_2^2.
 \ees
Gronwall inequality leads to the a priori estimate
 \be\label{AP}
\|(\nabla u,\nabla \om)\|_2^2\le \|(\nabla u_0,\nabla
\om_0)\|^2_2\exp\bigg\{C(p,\chi,\ga)\int^t_0\|u(s)\|_{p,\infty}^{2p/(p-3)}\md
s\bigg\}.
 \ee

In the case (b), we estimate $I_1$-$I_2$ in another way. Using the
generalized H\"older inequality and Young's inequality, we have
 \be\nonumber
I_1&\le& C\|\nabla u\|_{p,\infty}\|\nabla u\|_{\frac p{p-1},1}\le
C(p)\|\nabla u\|_{p,\infty}\|\nabla u\|^2_{\frac{2p}{p-1},2}\\
\nonumber &\le& C(p)\|\nabla u\|_{p,\infty}\|\nabla
u\|_2^{2-\frac3p}\|D^2 u\|_2^{\frac 3p}\\ \label{3.15} &\le&
\frac\chi{12}\|D^2u\|_2^2+C(p,\chi)\|\nabla
u\|_{p,\infty}^\frac{2p}{2p-3}\|\nabla u\|_2^2,
 \ee
where use has been made of the facts
 \bes
L^{\frac{2p}{p-1},2}(\mr^3)=(L^2,L^6)_{\frac{2p-3}{2p},2}(\mr^3)
 \ees
and
 \bes
\|\nabla u\|^2_{\frac{2p}{p-2},2}\le C\|\nabla
u\|_2^{2-\frac3p}\|D^2u\|_2^{\frac3p}.
 \ees
Arguing similarly, $I_2$ can also be estimated as follows
 \be\label{3.16}
I_2\le\frac\ga6\|D^2\om\|^2+C(p,\ga)\|\nabla
u\|_{p,\infty}^\frac{2p}{2p-3}\|\nabla \om\|_2^2.
 \ee
Inserting the estimates (\ref{3.15})-(\ref{3.16}) and (\ref{I3})
into (\ref{3.10}) and summing up $i$ from $1$ to $3$, and applying
Gronwall inequality, it reaches the a priori estimate
 \be\label{AP2}
\|(\nabla u,\nabla \om)\|_2^2\le \|(\nabla u_0,\nabla
\om_0)\|^2_2\exp\bigg\{C(p,\chi,\ga)\int^t_0\|\nabla
u(s)\|_{p,\infty}^{2p/(2p-3)}\md s\bigg\}.
 \ee
The above estimates are also valid for $p=\infty$ provided we modify
them accordingly.

 Combining the a priori estimates (\ref{AP}) and
(\ref{AP2}) with the energy inequality (\ref{EI}) and by standard
arguments of continuation of local solutions, we conclude that the
solutions $(u(t,x),\om(t,x))$ can be extended beyond $t=T$
provided that $u(t,x)\in L^q(0,T;L^{p,\infty}(\mr^3))$ for
$\frac2q+\frac3p\le 1$ with $3<p\le\infty$, or  $\nabla u(t,x)\in
L^q((0,T);L^{p,\infty}(\mr^3))$  for  $\frac2q+\frac3p\le 2$ with
$\frac32<p\le\infty$. The proof of Theorem \ref{thm2.1} is thus
competed.

 {\bf Proof of Theorem \ref{thm2.2}}: To prove the theorem we need
 the $L^4$ a priori estimate. For this purpose, we take the inner
 product of the first equation of (\ref{MMP}) with $|u|^2u$ and
 integrate by parts, it can be deduced that
 \be\nonumber
&&\frac14\frac{\md}{\md t}\|u\|_4^4+(\mu+\chi)\int_{\mr^3}|\nabla
u|^2|u|^2\md x+\frac12(\mu+\chi)\int_{\mr^3}|\nabla|u|^2|^2\md x\\
\label{3.20} &\le& 2\int_{\mr^3}|P||u|^2|\nabla u|\md
x+3\chi\int_{\mr^3}|w||u|^2|\nabla u|\md x,
 \ee
where we used the following relations by the divergence free
condition div$u=0$:
 \bes
&&\int_{\mr^3}u\cdot\nabla u\cdot|u|^2u\md x=\frac12\int_{\mr^3}u\cdot\nabla |u|^4\md x=0, \\
&&\int_{\mr^3}\Delta u\cdot|u|^2u\md x=-\int_{\mr^3}|\nabla
u|^2|u|^2\md x-\frac12\int_{\mr^3}|\nabla|u|^2|^2\md x,\\
&&\int_{\mr^3}\nabla\times\om\cdot|u|^2u\md
x=-\int_{\mr^3}|u|^2\om\cdot\nabla \times u\md
x-\int_{\mr^3}\om\cdot\nabla|u|^2\times u\md x,
 \ees
and
 \bes
|\nabla\times u|\le |\nabla u|, \ |\nabla|u||\le |\nabla u|.
 \ees
Using an argument similar to that used in deriving the estimate
(\ref{3.20}) it can be obtained for the second equation of
(\ref{MMP}) that
 \be\nonumber
&&\frac14\frac{\md}{\md
t}\|\om\|^4_4+\ga\int_{\mr^3}|\nabla\om|^2|\om|^2\md
x+\frac12\ga\int_{\mr^3}|\nabla|\om|^2|^2\md
x+\kappa\int_{\mr^3}|div\om|^2\md x+2\chi\int_{\mr^3}|\om|^4\md
x\\ \label{3.21}&& \le 3\chi\int_{\mr^3}|u||\om|^2|\nabla\om|\md
x.
 \ee
Combining estimates (\ref{3.20}) and (\ref{3.21}), we arrive at
 \be\nonumber
&&\frac14\frac{\md}{\md
t}(\|u||_4^4+\|\om\|^4_4)+(\mu+\chi)\int_{\mr^3}|\nabla
u|^2|u|^2\md x+\frac12(\mu+\chi)\int_{\mr^3}|\nabla|u|^2|^2\md x\\
\nonumber &&+\ga\int_{\mr^3}|\nabla\om|^2|\om|^2\md
x+\frac12\ga\int_{\mr^3}|\nabla|\om|^2|^2\md
x+\kappa\int_{\mr^3}|div\om|^2\md x+2\chi\int_{\mr^3}|\om|^4\md x\\
\nonumber &\le& 2\int_{\mr^3}|P||u|^2|\nabla u|\md
x+3\chi\int_{\mr^3}|w||u|^2|\nabla u|\md
x+3\chi\int_{\mr^3}|u||\om|^2|\nabla\om|\md x\\ \label{3.22}&=&
I+II+III.
 \ee
Applying H\"older and Young's inequalities for $II$, it follows that
 \be\label{3.23}
II\le \frac12\chi\int_{\mr^3}|\nabla u|^2|u|^2\md
x+C(\chi)\Big(\|u\|_4^4+\|\om\|_4^4\Big).
 \ee
Arguing similarly to above it can be derived for $III$ that
 \be\label{3.235}
III\le \frac12\ga\int_{\mr^3}|\nabla \om|^2|\om|^2\md
x+C(\ga)\Big(\|u\|_4^4+\|\om\|_4^4\Big).
 \ee
Concerning the term $I$, by virtue of the generalized H\"older
inequality (\ref{holder}) we have
 \bes
I&\le&
C(p)\||P|^{1/2}\|_{2p,\infty}\||P|^{1/2}|u|\|_{\frac{2p}{p-1},2}\|u\nabla
u\|_2\\ &\le&C(p)\|P\|^{1/2}_{p,\infty}\|u\nabla
u\|_2\|P\|^{1/2}_{\frac{2p}{p-1},2}\|u\|_{\frac{4p}{p-1},4}.
 \ees
Applying the divergence operator div to the first equation of
(\ref{MMP}), one formally has
 \be\label{Pu}
P=\sum_{i,j=1}^3R_iR_ju_iu_j,
 \ee
where $R_j$ denotes the $j-$th Riesz operator. By the boundedness
of Riesz operator on Lorentz space $L^{p,q}(\mr^3)$ for $1<p\le
q<\infty$, and applying the generalized H\"older inequality
(\ref{holder}) again to obtain that
 \be
\|P\|_{\frac{2p}{p-1},2}\le C\|u\|^2_{\frac{4p}{p-1},4}.
 \ee
So the term $I$ can be estimated as
 \be\label{3.24}
I\le C(p)\|P\|^{1/2}_{p,\infty}\|u\nabla
u\|_2\|u\|^2_{\frac{4p}{p-1},4}.
 \ee
In view of the real interpolation
 \bes
L^{\frac{4p}{p-1},4}(\mr^3)=(L^4,L^{12})_{\frac{2p-3}{2p},4}(\mr^3)
 \ees
and the Sobolev inequality we have
 \be\label{3.25}
\|u\|^2_{\frac{4p}{p-1},4}\le
C\|u\|_4^{\frac{2p-3}p}\|u\|_{12}^{\frac3p}\le
C\|u\|_4^{\frac{2p-3}p}\|u\nabla|u|\|_2^{\frac3{2p}}.
 \ee
Inserting the estimate (\ref{3.25}) into the estimate (\ref{3.24})
of $I$, one estimates $I$ as
 \be\label{3.26}
I\le C(p)\|P\|^{1/2}_{p,\infty}\|u\|_4^{\frac{2p-3}p}\|u\nabla
u\|_2^{\frac{2p+3}{2p}}\le
C(p,\chi)\|P\|_{p,\infty}^{\frac{2p}{2p-3}}\|u\|^4_4+\frac12\chi\|u\nabla
u\|^2_2.
 \ee
Inserting the estimates (\ref{3.26}) and (\ref{3.23})-(\ref{3.235})
of $I$, $II$ and $III$ into (\ref{3.22}) it follows that
 \bes
\frac14\frac{\md}{\md t}(\|u\|_4^4+\|\om\|^4_4)\le
C(p,\chi)\|P\|_{p,\infty}^{\frac{2p}{2p-3}}\|u\|^4_4+C(\chi,\ga)(\|u\|_4^4+\|\om\|^4_4).
 \ees
Gronwall's inequality implies that
 \be\label{estimate1}
\|u\|_4^4+\|\om\|^4_4\le
(\|u_0\|_4^4+\|\om_0\|^4_4)\exp\bigg\{C(p,\chi,\ga)\int^t_0(1+\|P\|_{p,\infty}^{\frac{2p}{2p-3}})\md
s\bigg\}.
 \ee

In the case (2), we estimate $I$ by another method. First $I$ also
equals to $I=\int_{\mr^3}\nabla P\cdot u|u|^2\md x$, so
 \bes
I&\le& \int_{\mr^3}|\nabla P|^{1/2}|\nabla P|^{1/2}|u|^3\md x\\
&\le& C\||\nabla P|^{1/2}\|_{2p,\infty}\||\nabla
P|^{1/2}|u|^3\|_{\frac{2p}{2p-1},1}\\ &\le& C(p)\|\nabla
P\|^{1/2}_{p,\infty}\|\nabla P\|^{1/2}_2\|u\|^3_{\frac
{12p}{3p-2},4},
 \ees
where we used the generalized H\"older inequality (\ref{holder}).
Noting the relation (\ref{Pu}) between $P$ and $u$, Riesz operator's
boundedness on $L^p$ for $1<p<\infty$ and the real interpolation
 \bes
L^{\frac{12p}{3p-2},4}(\mr^3)=(L^4,L^{12})_{\frac{p-1}{p},4}(\mr^3),
 \ees
The term $I$ can be estimated as
 \be\nonumber
I&\le& C(p)\|\nabla P\|^{1/2}_{p,\infty}\|u\nabla
u\|_2^{1/2}\|u\|_4^{\frac{3(p-1)}p}\|u\|_{12}^{\frac3p}\\&\le&\nonumber
C(p)\|\nabla P\|^{1/2}_{p,\infty}\|u\nabla
u\|_2^{\frac{p+3}{2p}}\|u\|_4^{\frac{3(p-1)}p}\\&\le&\label{3.28}
C(p,\chi)\|\nabla P\|^{\frac{2p}{3(p-1)}}_{p,\infty}\|u\|_4^4+
\frac12\chi\|u\nabla u\|_2^2,
 \ee
where the Sobolev embedding $L^6(\mr^3)\hookrightarrow \dot
H^1(\mr^3)$ and Young's inequality were used. Inserting the
estimates (\ref{3.28}) and (\ref{3.23})-(\ref{3.235}) of $I$, $II$
and $III$ into (\ref{3.22}) one also has
 \bes
\frac14\frac{\md}{\md t}(\|u\|_4^4+\|\om\|^4_4)\le
C(p,\chi)\|\nabla
P\|^{\frac{2p}{3(p-1)}}_{p,\infty}\|u\|^4_4+C(\chi,\ga)(\|u\|_4^4+\|\om\|^4_4).
 \ees
Applying Gronwall's inequality also to arrive at
 \be\label{estimate2}
\|u\|_4^4+\|\om\|^4_4\le
(\|u_0\|_4^4+\|\om_0\|^4_4)\exp\bigg\{C(p,\chi,\ga)\int^t_0(1+\|\nabla
P\|^{\frac{2p}{3(p-1)}}_{p,\infty})\md s\bigg\}.
 \ee
The above estimates are also valid for $p=\infty$ provided we modify
them accordingly.

Combining the a priori estimates (\ref{estimate1}) and
(\ref{estimate2}) with the energy inequality (\ref{EI}) and by
standard arguments of continuation of local solutions, we conclude
that the solutions $(u(t,x),\om(t,x))$ can be extended beyond
$t=T$ provided that $P(t,x)\in L^q((0,T);L^{p,\infty}(\mr^3))$,
for $\frac2q+\frac3p\le 2$ with  $\frac32<p\le\infty$, or $\nabla
P(t,x)\in L^q((0,T);L^{p,\infty}(\mr^3))$, for $\frac2q+\frac3p\le
3$ with  $1<p\le\infty$. We thus complete the proof of Theorem
\ref{thm2.2}.

 \vspace{0.4cm}

\textbf{Acknowledgements}  Author was partially supported by the
National Natural Science Foundation of China (No. 10771052), Program
for Science\&Technology Innovation Talents in Universities of Henan
Province (No. 2009HASTIT007), Doctor Fund of Henan Polytechnic
University (No.B2008-62).
 \vspace{0.4cm}


\end{document}